\documentclass[12pt,leqno]{article}
\usepackage[pdftex]{graphicx}
\usepackage{amsfonts,amssymb,amsmath}
\usepackage{ccaption,mathrsfs,amsmath,amssymb,amsthm}
\usepackage{amsmath, amscd, amsfonts,  tikz, times}

\newcounter{conjecture}
\newcounter{remark}
\newcommand{\eqnsection}{
    \renewcommand{\theequation}{\thesection.\arabic{equation}}
    \makeatletter
    \csname @addtoreset\endcsname{equation}{section}
    \makeatother}
\newtheorem{theorem}{Theorem}

\newtheorem{lemma}{Lemma}

\newcommand{\dd}{\delta}

\newcommand{\lar}{\longrightarrow}

\newcommand{\CC}{\mathbb{C}}

\def \be{\begin{equation}}
\def \ee{\end{equation}}
\def \bt{\begin{theorem}}
\def \et{\end{theorem}}
\def \bea{\begin{eqnarray}}
\def \eea{\end{eqnarray}}
\def \bas{\begin{eqnarray*}}
\def \eas{\end{eqnarray*}}

\def \ga{\gamma}

\def \th{\theta}

\def \ff{\infty}

\def \CCC{{\cal C}}

\def \({\left(}
\def \){\right)}

\def \bc{\begin{center} }
\def \ec{\end{center} }
\def \bs{\begin{slide} }
\def \es{\end{slide} }

\def\square{{\vcenter{\vbox{\hrule height.3pt
         \hbox{\vrule width.3pt height5pt \kern5pt
            \vrule width.3pt}
         \hrule height.3pt}}}}
\def\qed{{\hfill $\Box$ \bigskip}}

\newcounter{cccases}
\eqnsection
\begin{document}

\title{A probabilistic proof of the Open Mapping Theorem for analytic functions.}

\author{
\begin{tabular}{c}
\textit{Greg Markowsky} \\
Monash University \\
Victoria 3800, Australia \\
gmarkowsky@gmail.com
\end{tabular}}

\bibliographystyle{amsplain}

\maketitle \eqnsection \setlength{\unitlength}{2mm}

\begin{abstract}
\noindent The conformal invariance of Brownian motion is used to give a short proof of the Open Mapping Theorem for analytic functions.

%\vski

%2010 Mathematics subject classification: 30A99; 60J65.

%\vski

%Keywords: Analytic functions; Open Mapping Theorem; Brownian motion.
\end{abstract}

An important staple of the standard complex analysis curriculum is the {\it Open Mapping Theorem}, which is as follows.

\begin{theorem} \label{OMT}
Let $f$ be a nonconstant analytic function on an open set $W \subseteq \CC$. Then $f(W)$ is open.
\end{theorem}
The standard proof, contained in virtually any complex analysis textbook, employs contour integration and the argument principle (or, equivalently, Rouch\'e's Theorem). In this note, we give a quick and intuitive proof of the theorem which eschews contour integration, utilizing instead the conformal invariance of Brownian motion and some simple topology. The following is what is referred to as conformal invariance, although it might also be referred to as {\it analytic invariance}, since injectivity is not required.

\begin{theorem} \label{holinv}
Let $f$ be analytic and nonconstant on a domain $U$, and let $a \in U$. Let $B_t$ be a planar Brownian motion in $U$ starting at $a$ and stopped at a stopping time $\tau$. Then there is a time change $t \lar C_t$ such that $\hat B_t = f(B_{C_t})$ is Brownian motion starting at $f(a)$ and stopped at the stopping time $C^{-1}_\tau$, where $C^{-1}_s = \inf\{t \geq 0: C_t \geq s\}$.
\end{theorem}

See \cite[Thm. V.2.5]{revyor} or \cite[Sec. 2.12]{durBM} for a proof of this theorem, which is originally due to Paul L\'evy. For the proof of Theorem \ref{OMT}, we need the following lemma.

\begin{lemma} \label{hitexit}
Let $a \in \CC$, and for any $\dd > 0$ let $\tau_\dd = \inf \Big\{t \geq 0 : B_t \in \{|z - a| = \dd\}\Big\}$. Fix $r>0$, and let $V$ be an open set contained in $D(a,r)$, where $D(a,r) = \{|z-a| < r\}$. Let $B_t$ be a Brownian motion starting at $a$. Then $P(B_t \in V \mbox{ for some } t \leq \tau_r) > 0$.
\end{lemma}

{\bf Proof:} Since $V$ is open, it contains a circular arc $\CCC = \{a + r'e^{i\th}: \th_1 < t < \th_2\}$ with $0<r' < r$ and $0 < \th_2 - \th_1 \leq 2 \pi$. We have

\vspace{-.2in}

\begin{equation} \label{}
P(B_t \in V \mbox{ for some } t \leq \tau_r) \geq P(B_t \in \CCC \mbox{ for some } t \leq \tau_r) \geq P(B_{\tau_{r'}} \in \CCC) = \frac{\th_2 - \th_1}{2\pi};
\end{equation}
note that the last equality holds by the rotation invariance of $B_t$. \qed

{\bf Proof of Theorem \ref{OMT}:} Let $v \in f(W)$, and choose $a \in W$ such that $f(a) = v$. By the Identity Theorem for analytic functions (see \cite[Thm. 10.18]{rud}), we can choose a small $r>0$ such that $f(w) \neq v$ on $\{|z-a|=r\}$ and $\bar D(a,r) \subseteq W$, where $\bar D(a,r) = \{|z-a| \leq r\}$. If we start a Brownian motion at $a$ and stop it at the hitting time $\tau$ of $\{|z-a|=r\}$, then $\hat B_t = f(B_{C_t})$ is a Brownian motion starting at $v$ and stopped at the stopping time $C^{-1}_\tau$, which satisfies $f(\hat B_{C^{-1}_\tau}) = f(B_{\tau}) \in f(\{|z-a|=r\})$. Let $\ga$ denote the curve $f(\{|z-a|=r\})$. In traveling from $v$ to $\ga$, $\hat B_t$ must cross first the circle $\{|w-v|=m\}$, where $m = \inf \{|w-v|: w \in \ga \} > 0$; so if we let $\hat \tau$ be the first time $\hat B_t$ hits $\{|w-v|=m\}$ then $C^{-1}_\tau \geq \hat \tau$ a.s.

\vspace{-1.5in}

\includegraphics[width=180mm,height=145mm]{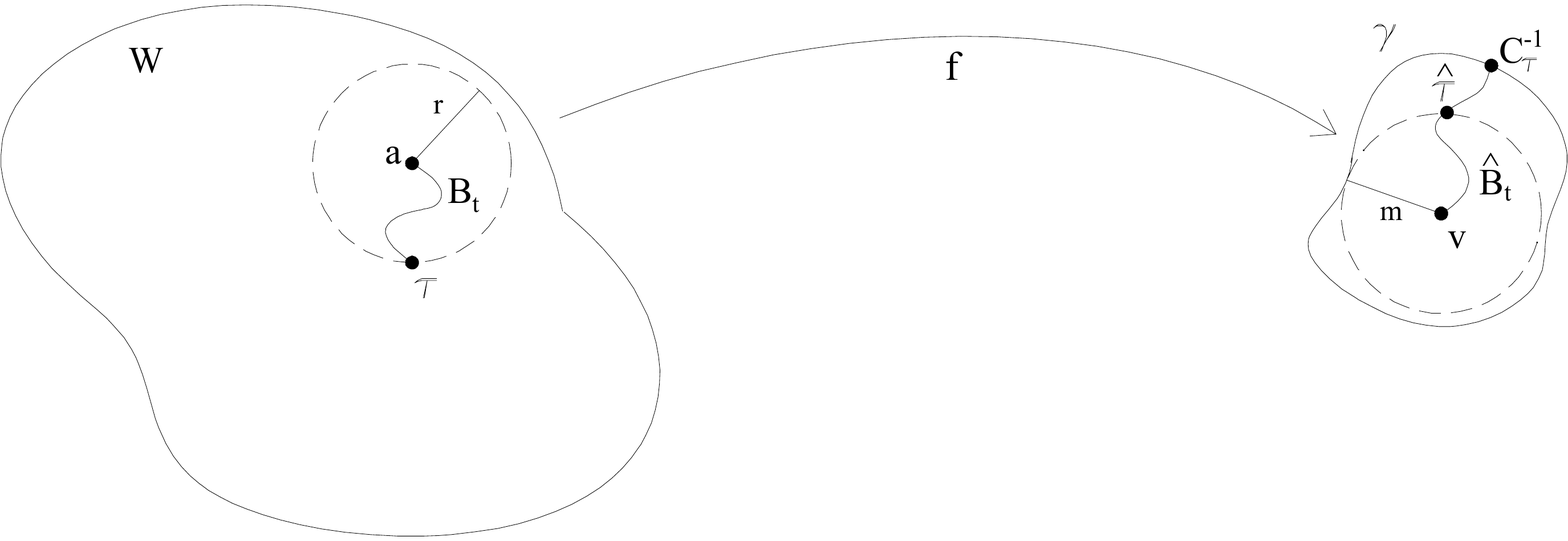}

\vspace{-1.55in}

Lemma \ref{hitexit} now shows that $\hat B_t = f(B_{C_t})$ hits any open set in $\{|w-v|<m\}$ with positive probability before time $\hat \tau$, and this shows in particular that $\{|w-v|<m\}$ is contained in the closure of $f(\bar D(a,r))$. However, $\bar D(a,r)$ is a compact set, so its image under the continuous map $f$ is compact and thus closed, and therefore $\{|w-v|<m\} \subseteq f(\bar D(a,r))\subseteq f(W)$. This proves that $f(W)$ is open. \qed

{\bf Remark:} In the picture above, $\ga$ appears as a Jordan curve encircling $v$. Of course, if $f$ is not injective in $\{|z-a| \leq r\}$ then the curve may have self-intersections; a larger point is that the proof above does not require the fact that $\ga$ separates $v$ from $\ff$, but uses only that the curve does not intersect $\{|w-v|<m\}$. It is also not true in general that $C^{-1}_\tau$ must be the hitting time of $\ga$, as the picture suggests, but again the proof above does not require this, using only the fact that $f(\hat B_{C^{-1}_\tau}) \in \ga$.

\section{Acknowledgements}

I'd like to thank Andrea Collevecchio for helpful discussions.

\bibliographystyle{alpha}
\bibliography{CABMbib}

\end{document}